\documentclass[12pt]{amsart}

\newcounter{defcounter}
\usepackage{amssymb,amsmath,amscd,graphicx, enumitem,
latexsym,amsthm}
\usepackage{amssymb,latexsym,eufrak,amsmath,amscd,graphicx}
  \usepackage[all]{xy}
  \usepackage{pdfsync}
\theoremstyle{plain}

\newtheorem*{theoremnn}{Theorem}
\newtheorem{proposition.definition}[theorem]{Proposition/Definition}
\newtheorem*{conjecturenn}{Conjecture}
\newtheorem*{exercisenn}{Exercise}

\theoremstyle{definition}

\newcommand{\noi}{\noindent}
\newcommand{\PP}{\mathbf{P}}

\newcommand{\ZZ}{\mathbf{Z}}
\newcommand{\CC}{\mathbf{C}}
\newcommand{\QQ}{\mathbf{Q}}

\newcommand{\OO}{\mathcal{O}}

\newcommand{\fra}{\mathfrak{a}}

\newcommand{\frq}{\frak{q}}

\newcommand{\HH}[3]{H^{{#1}} \big( {#2} , {#3}
\big) }

\newcommand{\lct}{\textnormal{lct}}

\newcommand{\mult}{\textnormal{mult}}

\newcommand{\Linser}[1]{| \mspace{1.5mu} {#1}
\mspace{1.5mu} |}
\newcommand{\linser}[1]{\Linser{  {#1}  }}

\newcommand{\charr}{\textnormal{char}}

\numberwithin{theorem}{section}

\begin{document}

\title
{Some remarks on the work of Lawrence Ein}

 \author{Robert Lazarsfeld}
  \address{Department of Mathematics, Stony Brook University, Stony Brook, New York 11794}
 \email{{\tt robert.lazarsfeld@stonybrook.edu}}
 \thanks{Research of the  author partially supported by NSF grant DMS-1439285.}

\maketitle

 \subsection*{Introduction}

To set the stage for the conference, I gave an informal survey  -- seasoned with  stories and reminiscenses -- of Ein's mathematical work. The present pages constitute a write-up of these remarks.  I'm not sure  how well such a talk will translate into print, but I appreciate this opportunity to express my esteem for Lawrence's mathematics as well as my gratitude to him for our long collaboration.

\subsection*{Vector Bundles}  
  Lawrence arrived in Berkeley in 1977, where he joined an active group of students around Hartshorne and Ogus that included Ziv Ran, Mei-Chu Chang and Tim Sauer. At the time Hartshorne in particular was very interested in the geometry of vector bundles on projective space, and Lawrence's first papers were devoted to questions in this direction. In \cite{Ein1} -- which became his thesis -- Ein extended to characteristic $p > 0$ some of the classical results of Van de Ven and others about  bundles on $\PP^n(\CC)$. For those of us who think of Lawrence as mainly traveling in  characteristic zero, it is interesting that this first paper involves some serious positive characteristic geometry. For example, Ein classifies uniform $n$-bundles on $\PP^n(k)$, and finds that when $\charr \ k > 0$ some Frobenius pullbacks show up.  In another paper \cite{Ein2},  Ein, Hartshorne and Vogelaar prove analogues for rank three bundles of classical results of Barth on restruictions of rank $2$ steble bundles to hyperplanes.

Starting in the mid 1980's, Lawrence shifted towards \textit{applying} ideas involving vector bundles to study concrete geometric problems. This would become characteristic of some of his most interesting work during the 1980s. For example, in \cite{Ein3}, Ein  established a Noether-type theorem for the Picard groups of surfaces arising as the degeneracy locus associated to general sections of a suitably positive vector bundle. In the same paper he showed that the Le Potier vanishing theorem leads to a very quick proof of a result of Evans and Griffith \cite{Evans.Griffith} that if $E$ is a vector bundle of rank $e$ on a complex projective space $\PP$ that satisfies the vanishing
\[   \HH{i}{\PP}{E(k)} \ = \ 0 \ \ \text{ for all } \ 0 < i < e \ , \ k \in \ZZ, \]
then $E$ is a direct sum of line bundles. Evans and Griffith had deduced this as a consequence of their deep  algebraic results on syzygy modules, and it was very nice to have a quick geometric proof of the statement. 

\subsection*{Subvarieties of general hypersurfaces and varieties with small duals}  Lawrence's most influential work during the period 1985--1990 are arguably his papers \cite{Ein4}, \cite{Ein5} on varieties with small duals, and his results \cite{Ein6}, \cite{Ein7} on subvarieties of very general complete intersections.

Let $X \subseteq \PP^r$ be a smooth complex projective variety of dimension $n$. Recall that the \textit{dual variety} of $X$ is the set of hyperplanes tangent to $X$ at some point:
\[   X^* \ =_{\text{def}} \ \big \{ H \subseteq \PP^r \mid H \text{ is tangent to } X \big \} \ \subseteq \ \PP^{r^*}.\]
One expects $X^*$ to be a hypersurface in $\PP^{r*}$, but sometimes it has smaller dimension. Zak had established using the Fulton-Hansen connectedness theorem that in any event $\dim X^* \ge \dim X$. Zak's work brought renewed attention to the classical question of trying to say something about those smooth varieties whose duals are exceptionally small.

In his first paper \cite{Ein4}, Lawrence shows that the only smooth varieties $X \subseteq \PP^r$  with
\[   \dim X \ = \ \dim X^* \ \le \ \tfrac{2}{3} r \]
are four classically known examples.\footnote{Hartshorne's famous conjecture on complete intersection predicts that if $\dim X > \tfrac{2}{3} r$ then $X$ should be a complete intersection, and in particular will have a non-degenerate dual.} In his second paper \cite{Ein5}, he classifies all varieties of dimension $<6$ whose duals are degenerate. 
Ein starts by fixing a general tangent hyperplane $H \subseteq \PP^r$ to $X$. Then (as was classically understood) the \text{contact locus}
\[   L \ =_{\text{def}}\ \big \{ x \in X \mid \text{ $H$ is tangent to $X$ at $x$ } \big \}\]
is a \textit{linear} space of dimension $k = r - 1 - \dim X^*$. Lawrence's very nice idea is to study the normal bundle
\[   N \  = \ N_{L/X} \]
to $L$ in $X$.  This is a  bundle of rank $n - k$ on $L = \PP^k$, and Ein proves that
\[    N \ \cong \ N^* \otimes \OO_L(1),  \tag{*}
\]
which he then combines with input from the geometry of vector bundles on projective space. Note that (*) already implies Landman's theorem that if $k \ge 1$, then
\[ k \ \equiv \ \dim X \pmod{2}. \]
Some related ideas appear in Ein's paper \cite{ESB} with Shepherd-Barron concerning special Cremona transformations.

In 1986, Herb Clemens \cite{Clemens} proved a lower bound on the geometric genus of a curve on a very general hypersurface in projective space. Ein's two very influential papers \cite{Ein6} and \cite{Ein7} gave a large generalization of this result.
\begin{theoremnn}
Let $X \subseteq \PP^n$ be a very general complete intersection of type $(m_1, \ldots, m_e)$, and assume that 
\[ \sum m_i \ \ge \ 2n + 1 -e \]
Then any subvariety of $X$ is of general type.
\end{theoremnn}
\noi Very roughly speaking, Lawrence's  idea was that the PGL action on projective space allows one to produce sections of a suitable twist of the normal or canonical bundle of such a subvariety. For example, he shows that if $Z \subset X$ is a smooth subvariety that deforms with $X$, and if $N= N_{Z/X}$ denotes the normal bundle of $Z$ in $X$,  then $N(1)$ is
globally generated. By adjunction, it follows that $X$ must be of general type in a suitable range of degrees. Voisin \cite{V1} subsequently gave a particularly clean formulation of the argument. 

\subsection*{Linear series on higher dimensional varieties} I first got to know Lawrence during the 1981-82 academic year, when we were both  at the IAS. 
We kept in close mathematical touch after that, but we only started actively collaborating in the late 1980's. The initial fruits of these efforts were the three papers \cite{BEL}, \cite{SAD} and \cite{PALS} (the first with Bertram) dealing with linear series on higher dimensional varieties. Ein and I had both been interested over the years in questions about linear series, and this was the time at which higher-dimensional geometry was really beginning to flower. So it seemed reasonable to see whether one could say something in higher dimensions about issues that had attracted attention for curves and surfaces, notably syzygies and Castelnuovo-Mumford regularity.  The idea of \cite{BEL} was to use vanishing theorems to study  the regularity of varieties defined by equations of given degrees, while \cite{SAD} proved analogues for all smooth varieties of results of Mark Green on syzygies of curves and Veronese varieties.  This work was completed around 1990, and we were happy that it at least let us get our toes in the water of higher dimensions. 

The paper \cite{PALS}, dealing with global generation of linear series on threefolds, gave us our first practice in using the cohomological techniques of Kawamata--Reid--Shokurov involving the Kawamata--Viehweg vanishing theorem for $\QQ$-divisors. We came to this work through a fortuitous sequence  of  events, illustrating in an amusing way the role that chance  sometimes plays in setting the direction of one's research. It's perhaps worth telling the story.

In the late 1980's, Fujtia proposed a far-reaching generalization of the classical fact that  a line bundle of degree $\ge 2g$ (or $\ge 2g +1$) on a curve of genus $g$ is globally generated (or very ample):
\begin{conjecturenn} Let $A$ be an ample line bundle on a smooth projective variety $X$ of dimension $n$. Then:
\begin{itemize}
\item[$(i)$.] $K_X + (n+1)A$ is basepoint-free;
\vskip 5pt
\item[$(ii)$.] $K_X +(n+2)A$ is very ample.
\end{itemize}
\end{conjecturenn}
\noi Fujita was also influenced by then recent work of Igor Reider \cite{Reider}, who had studied adjoint line bundles $K_X + L$ on surfaces, and whose results implied the $\dim X = 2$ case of the conjecture. The first breakthrough in Fujita's conjecture occured around 1990, when Demailly \cite{Demailly} applied $L^2$-methods to prove:
\begin{theoremnn}[Demailly] In the situation of Fujita's conjecture, 
\[   2K_X + 12n^nA \text{ is very ample}.\]
\end{theoremnn} 
\noi While the numerics were rather far from what one expected the result created a sensation, going as it did far beyond what algebraic methods gave at the time. Moreover, Demailly's work inaugurated a very   fruitful interaction between analytic and algebraic geometry that continues to this day.

In 1991, Lawrence and I attended a meeting at Oberwolfach at which Siu gave a report aimed at algebraic geometrers on Demailly's work and some extensions thereof. (His abstract is reproduced in Figure 1.) 
During the course of his talk, Siu posed as a challenge the
\begin{exercisenn}
Let $L$ be an ample line bundle on a smooth projective $n$-fold $X$, and   fix a point $x \in X$. Suppose that there exists $k \gg 0$ and  a divisor
 $ D \in \linser{kL}$
 such that
\[  \mult_x(D) \ge nk, \]
while $\mult_y(D) < k$ for $y$ in a punctured neighborhood of $x$. Then $K_X + L$ is free at $x$.

\end{exercisenn}
 \begin{figure}  \label{Siu}
 \includegraphics[scale = .50]{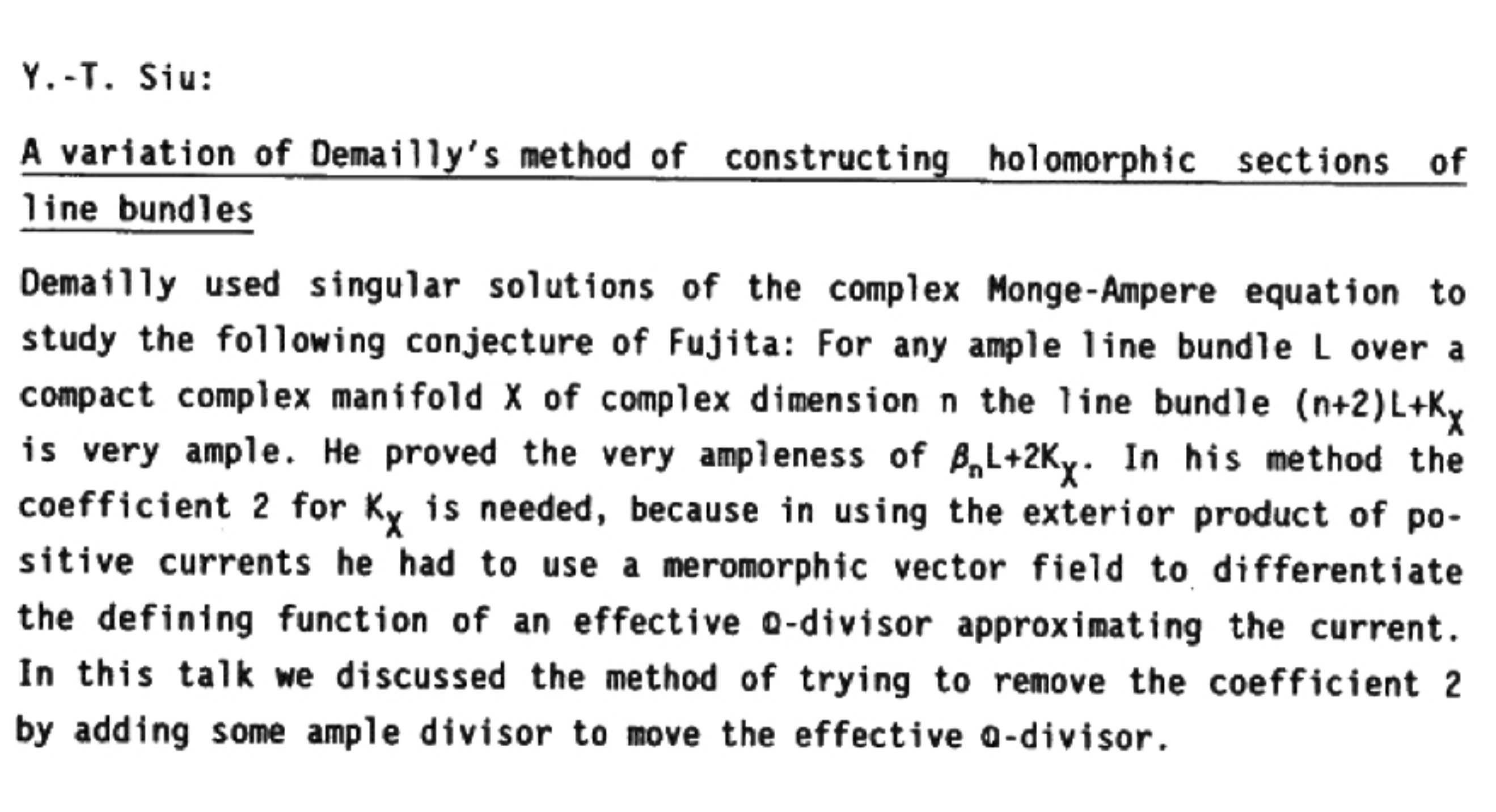}  
 \caption{Siu's 1991 Oberwolfach Abstract}  \end{figure}

 \noi As soon as Siu put the statement on the board, H\'el\`ene Esnault and Eckart Viehweg, who were in the audience, realized that one should be able to solve it using vanishng theorems for $\QQ$-divisors, along the lines of their paper \cite{EV} from some years back. By the end of the hour they had worked out the details, which they kindly explained to Lawrence and me over the following days. (At that time, the Kawamata--Viehweg vanishing theorem still seemed rather exotic, and its geometric content was not well understood.)
 
 Lawrence and I noticed that a parameter count suggested that one might use Siu's exercise to reprove Reider's results on surfaces, and somewhat later -- after we learned more of cohomological arguments of Kawamata--Reid--Shokurov -- we got this working. But we assumed that the experts were  aware of this approach,  which discouraged us from trying to tackle higher dimensions. A couple of months  later, however, I was chatting with J\'anos Koll\'ar at a UCLA--Utah workshop, and when I mentioned to him in passing this  well-known argument, J\'anos seemed rather surprised to hear about it.\footnote{In fact he asked me to show him the proof, but  it wasn't fresh on my mind and that part of the conversation didn't go very well...} At that point  Lawrence and I  realized  that people actually hadn't tried to use Kawamata's ``X-method" to give effective results, and after a couple of months of very intense work we were able to prove the global generation statement in dimension three. With the hindsight of subsequent developments such as \cite{AS}, \cite{Kol}, \cite{Kawamata},  the resulting paper \cite{PALS} today looks hopelessly   overcomplicated and dated, but we like to think it played a useful role by providing   proof of concept for extracting reasonably sharp effective statements from  the cohomological apparatus of the minimal model program. 
 
 \subsection*{Multiplier ideals} Work of Demailly, Siu and others (eg \cite{Demailly}, \cite{AS}, \cite{Siu1}, \cite{Siu2}) gave dramatic evidence of the power of multiplier ideals and the vanishing theorems they satisfy. Starting in the mid 1990's, Lawrence and I made a systematic attempt to look for concrete applications of this machinery. A first one appeared in \cite{STIV}, where we combined generic vanishing theorems with some  cousins of multiplier ideals to prove a conjecture about the singularities of theta divisors.  Another application appears in \cite{GEN}, where we applied   an algebro-geometric analogue of Skoda's theorem to establish variants of Koll\'ar's results on the effective Nullstellensatz. But my favorite result from this period is the paper \cite{UBSP} with Karen Smith, which used multiplier ideals to prove a rather surprising comparison between the symbolic and ordinary powers of a radical ideal on a smooth variety.
 
 One of the things that Ein teaches his collaborators is that it can be a useful exercise to try to use new methods to reprove old results. During Spring 2000, Lawrence visited Michigan while on sabbatical, and with Karen  we decided to look for applications of the circle of ideas around the subadditivity theorem \cite{DEL}. We came across an old result of Izumi \cite{Izumi} to the effect that ideals associated to a divisorial valuation are contained in growing powers of the ideal of its support, and saw that this followed immediately from subadditivity. But in fact, the argument seemed to  prove  more: namely, it showed that if $\frq \subseteq \OO_X$ is the ideal sheaf of a reduced subvariety $Z \subseteq X$ of a smooth variety of dimension $d$, then 
 \[   \frq^{(md)} \ \subseteq \ \frq^m, \] 
 where the symbolic power on the left denotes the germs of functions having multiplicity $\ge md$ at a general point of $Z$.  At first  it seemed unclear whether one could really expect such a statement to be true, but luckily the proof was only a couple of lines long and so there weren't many places for an error to hide. Soon thereafter, Hochster and Huneke \cite{HH} gave a different approach to this comparison via tight closure.\footnote{The ideas surrounding Izumi's theorem were further developed in \cite{AbhVals}.}

 \subsection*{Singularities}  During the first decade of the 2000's, Lawrence's focus shifted to the geometry of the singularities that arise in birational geometry. Working with Musta\c t\u a, de Fernex and others, he resolved a number of long-standing problems in the area. For example  in \cite{EMY}  Ein, Musta\c t\u a and Yasuda used the geometry of jets to prove the conjecture on  inversion of adjunction on smooth varieties. At the end of this period, inspired by some ideas of Koll\'ar, Ein, de Fernex and Musta\c t\u a  \cite{dFEM2} proved the celebrated ACC conjecture on smooth varieties:
 \begin{theoremnn}
 The collection of rational numbers that can occur as log-canonical thresholds of divisors on a smooth $n$-fold does not contain any cluster points from the left.
 \end{theoremnn}
 \noi This was later extended to varieties with mild singularities by Hacon, McKernan and Xu \cite{HMX}. 
 
 Some other papers during these years established interesting new relations among classical invariants, and applied them to questions of birational rigidity. For example in \cite{dFEM1}  Ein, de Fernex and Musta\c t\u a  prove
 \begin{theoremnn}
 Let $\fra \subseteq \OO_X$ be an ideal of finite colength on a smooth variety $X$ of dimension $n$. Then
  \[  e(\fra) \ \ge \ \left( \frac{n}{\textnormal{lct}(\fra)}\right)^n. \]
 \end{theoremnn}
 \noi Here $e(\fra)$ denotes the classical Samuel multiplicity of $\fra$, and $\lct(\fra)$ denotes the log-canonical threshold of a general element in $\fra$. The proof in  \cite{dFEM1} used a degeneration to reduce the question to the case of monomial ideals. Yuchen Liu \cite{Liu}  recently found a very different argument via the theory of K-stability. In \cite{dFEM1.5}, Ein, de Fernex and Musta\c t\u a  applied variants of their inequality to prove the birational rigidity of smooth hypersurfaces of degree $n$ in $\PP^n$ when $4 \le n \le 12$. De Fernex extended this result to all dimensions in \cite{dF}.

 \subsection*{Asymptotic syzygies} Over the last few years, Lawrence and I have returned to the circle of questions around the defining equations of projective varieties and the syzygies among them. Specifically, we have tried to understand what happens syzygetically  as $d \to \infty$ when one takes a fixed projective variety $X$ and embeds it by  line bundles of the form 
 \[ L_d \ = \ d A \, + \, P, \]
where $A$ is ample and $P$ is arbitrary. The paper \cite{ASAV}  shows in effect that essentially all the Koszul cohomology groups that could be non-zero  asymptotically become so. 
 It was important technically to deal not just with the Koszul cohomology groups $K_{p,q}(X; L_d)$ of $L_d$, but to work more generally with the groups $K_{p,q}(X, B; L_d)$ in which one allows a fixed twisting line bundle. Combining this perspective with Voisin's Hilbert schematic approach to syzygies \cite{Voisin1}, \cite{Voisin2} led to a remarkably simple proof in \cite{Gonality} of an  old conjecture with Mark Green from \cite{GL} concerning the syzygies of curves of large degree:
 \begin{theoremnn}
 One can read off the gonality of a curve $C$ from the grading of its resolution in any one embedding of sufficiently large degree. 
 \end{theoremnn}
 \noi Rathmann \cite{Rathmann} subsequently found a very nice argument giving an effective bound: it suffices in the Theorem to consider embeddings of degree $\ge 4 g(C) + 1$.

 \subsection*{Epilogue} This  overview has omitted many of Ein's contribitions, and we can be sure that  there will be  more to come in the future. Still, I hope these remarks have conveyed something of the breadth  of Lawrence's work and vision. Through his research and his generosity in sharing an encyclopedic  knowledge of the field, Ein has had an enormous and continuing impact on several generations of algebraic geometers. Working with Lawrence has certainly been the most important influence on my own mathematical career. It is a joy to have this occasion to wish him Happy Birthday!
 
 %
 %
 %
 %

 \end{document}